%
%
%
\magnification=1200
\input amstex
\documentstyle{amsppt}
\def\cfrac#1#2{\dfrac{\strut#1}{#2}\kern-\nulldelimiterspace}
\topmatter
\title
The Impact of Stieltjes' Work on Continued Fractions and Orthogonal Polynomials
\endtitle
\author
Walter Van Assche
\endauthor
\affil
Katholieke Universiteit Leuven
\endaffil
\address
Department of Mathematics, Katholieke Universiteit Leuven,
Celestijnenlaan 200\,B, B-3001 Leuven (Heverlee), BELGIUM
\endaddress
\email
fgaee03\@cc1.kuleuven.ac.be
\endemail
\thanks
The author is a Research Associate of the Belgian National Fund for Scientific
Research
\endthanks
\leftheadtext{Walter Van Assche}
\rightheadtext{Impact of Stieltjes' work}

\abstract
Stieltjes' work on continued fractions and the orthogonal polynomials related
to continued fraction expansions is summarized and an attempt is made to
describe the influence of Stieltjes' ideas and work in research done after his
death, with an emphasis on the theory of orthogonal polynomials.
\endabstract

\toc
\head {} Introduction \endhead
\head 1. Stieltjes Continued Fraction \endhead
\head 2. Moment Problems \endhead
   \subhead 2.1 The Stieltjes moment problem \endsubhead
   \subhead 2.2 Other moment problems \endsubhead
   \subhead 2.3 Recent extensions of the moment problem \endsubhead
\head 3. Electrostatic Interpretation of Zeros  \endhead
   \subhead 3.1 Jacobi polynomials \endsubhead
   \subhead 3.2 Laguerre and Hermite polynomials \endsubhead
   \subhead 3.3 Extensions  \endsubhead
   \subhead 3.4 Logarithmic potential theory \endsubhead
\head 4. Markov-Stieltjes Inequalities \endhead
\head 5. Special Polynomials \endhead
   \subhead 5.1 Legendre polynomials \endsubhead
   \subhead 5.2 Stieltjes polynomials \endsubhead
   \subhead 5.3 Stieltjes-Wigert polynomials \endsubhead
   \subhead 5.4 Orthogonal polynomials related to elliptic functions
\endsubhead \endtoc
\endtopmatter

\document
\head Introduction  \endhead
The memoir {\it Recherches sur les fractions continues}, published
posthumously in the Annales de la Facult\'e des Sciences de Toulouse --- a
journal of which Stieltjes was one of the first editors --- and a great
number of other papers by Stieltjes contain a wealth
of material that still has a great impact on contemporary research,
especially on the theory of orthogonal polynomials. The general
theory of orthogonal polynomials really started with the investigations
of Chebyshev and Stieltjes. The impact of the work of Chebyshev and his student
Markov has already been described by Krein \cite{53}. Here we give an attempt
to discuss some of Stieltjes' contributions and the impact on later work.
Orthogonal polynomials offer a variety of results and applications. The
bibliography \cite{91} up to 1940 consists of 1952 papers by 643 authors. Even
now interest in orthogonal polynomials is enormous. One of the reasons is that
orthogonal polynomials seem to appear in a great variety of applications. Their
use in the numerical approximation of integrals was already pointed out by
Gauss and further extended by Christoffel \cite{19} \cite{30} and
Stieltjes \cite{95}. The Pad\'e table \cite{73}
for the approximation of a function by rational functions is very closely
related to continued fractions and Stieltjes' work may be considered as
one of the first proofs of convergence in the Pad\'e table \cite{74}.
In 1954 Lederman and Reuter \cite{56} and in
1957 Karlin and McGregor \cite{49} showed that the transition probabilities
in a birth and death process could be expressed by means of a Stieltjes integral
of orthogonal polynomials. Even in pure mathematics there seems to be a natural
framework where orthogonal polynomials come into play: representations for
certain Lie groups very often are in terms of special functions, in particular
orthogonal polynomials (see e.g. Vilenkin \cite{117}). Recently this has also
been observed for quantum groups \cite{52}.
Discrete orthogonal polynomials have useful
applications in the design of association schemes and the proof of nonexistence
of perfect codes and orthogonal polynomials on the unit circle have a close
connection with digital signal processing. The proceedings of the NATO
Advanced Study Institute on ``Orthogonal Polynomials and their Applications''
(Columbus, Ohio 1989) \cite{69} gives excellent contributions to each of these
aspects of orthogonal polynomials and is strongly recommended.

Stieltjes' work has already been discussed by Cosserat \cite{22} shortly
after Stieltjes' death in 1894. In these notes we will try to estimate the
value of the investigations by Stieltjes a century later. Let me also
mention Brezinski's book on the history of continued fractions
\cite{10, Chapter 5, Section 5.2.4 on pp.\ 224--235} where Stieltjes'
work on continued fractions is shown in its historic context.

\head 1. Stieltjes Continued Fraction \endhead

The object of his main work \cite{105} is the study of the continued fraction
$$    \cfrac{1}{c_1z + \cfrac{1}{c_2 + \cfrac{1}{c_3z + \cdots +
\cfrac{1}{c_{2n} + \cfrac{1}{c_{2n+1}z + \cdots}}}}}, \tag 1.1 $$
which is nowadays known as a Stieltjes continued fraction or $S$-fraction.
Stieltjes only considers the case where $c_k > 0$ $(k=1,2,\ldots)$.
In general an $S$-fraction is any continued fraction of the form \thetag{1.1}
in which all $c_k$ are different from zero or any continued fraction which can
be obtained  from it by an equivalence transformation or change of variable
\cite{119, p.\ 200}. The $S$-fraction \thetag{1.1} can be transformed by
contraction to a $J$-fraction
$$  \cfrac{a_0^2}{z-b_0 - \cfrac{a_1^2}{z-b_1 - \cfrac{a_2^2}{z-b_2 - \cdots -
\cfrac{a_{n-1}^2}{z-b_{n-1} - \cfrac{a_n^2}{z - b_n - \cdots}}}}}, \tag 1.2 $$
with $a_0^2 = 1/c_1, b_0 = -1/(c_1c_2)$ and
$$  a_n^2 = \frac{1}{c_{2n-1}c_{2n}^2c_{2n+1}}, \quad
    b_n = - \frac{1}{c_{2n}c_{2n+1}} - \frac{1}{c_{2n+1}c_{2n+2}},
    \qquad k=1,2,\ldots    . $$
The positivity of the $c_k$, as imposed by Stieltjes, clearly puts some
constraints on the coefficients $a_k, b_k$ e.g., $a_k^2 > 0$ and $b_k < 0$.
A $J$-fraction can be regarded
as being generated by the sequence of transformations
$$   t_0(w) = \frac1w, \qquad t_k(w) = z - b_{k-1} - \frac{a_k^2}{w},\qquad
k=1,2,\ldots $$
The superposition $t_0(t_1(\cdots(t_n(w))\ldots))$ for $w=\infty$ is then the
$n$-th approximant or $n$-th convergent of the fraction \thetag{1.2}. This is a
rational function of the variable $z$  and we have
$$   t_0(t_1(\cdots(t_n(\infty))\ldots)) = \frac{1}{a_1} \,
\frac{p_{n-1}^{(1)}(z)}{p_n(z)}, $$
where both the denominator polynomials $p_n(z)$ $(n=0,1,2,\ldots)$ and
numerator polynomials $p_{n-1}^{(1)}(z)$ $(n=0,1,2,\ldots)$ are solutions
of the three-term recurrence relation
$$ zr_n(z) = a_{n+1}r_{n+1}(z) + b_n r_n(z) + a_nr_{n-1}(z), \qquad n \geq 0
                     \tag 1.3 $$
with initial condition
$$ p_{-1}(z) = 0,\ p_0(z) = 1, \qquad p_{-1}^{(1)}(z) = 0,\ p_0^{(1)}(z) = 1. $$
The convergents of the $S$-fraction are such that the $2n$-th convergent of
\thetag{1.1} is equal to the $n$-th convergent of \thetag{1.2}.
If the denominator $p_n(z)$ vanishes for at most a finite number of integers $n$
and if $\lim_{n \rightarrow \infty} p_{n-1}^{(1)}(z)/p_n(z) = f(z)$
exists, then the $J$-fraction converges to $f(z)$.
Stieltjes gave a general theory of $S$-fractions (and consequently of
$J$-fractions) with $c_k > 0$ $(k=1,2,\ldots)$,
dealing with questions of convergence and he showed a close
connection with asymptotic series in terms of a given sequence of moments (see
also the next section).

One of the most important facts in the theory is that the denominators
$p_n(-x)$ $(n=0,1,2,\ldots)$ form a sequence of orthonormal polynomials
on $[0,\infty)$ i.e.,
there is a positive measure $\mu$ on $[0,\infty)$ such that
$$  \int_0^\infty p_n(-x)p_m(-x) \, d\mu(x) = \delta_{m,n} . $$
The support of the measure $\mu$ is in $[0,\infty)$ precisely because
Stieltjes assumes the coefficients $c_k$ of the $S$-fraction \thetag{1.1} to
be positive.
Stieltjes showed that such orthogonal polynomials have zeros with interesting
properties. He proved that all the zeros of $p_n(-x)$ are real, positive and
simple; moreover
the zeros of $p_n(-x)$ interlace with the zeros of $p_{n-1}(-x)$ but also with
the zeros of $p_{n-1}^{(1)}(-x)$. The latter property shows that the
convergent $p_{n-1}^{(1)}(z)/p_n(z)$ is a rational function with $n$ real
and negative
poles and positive residues. These properties are now quite classical and
of great use for numerical quadrature.  The property of
orthogonality is crucial in these considerations (but Stieltjes never uses
this terminology). A famous and very important result in the theory of
orthogonal polynomials on the real line is the following result:
\proclaim{Theorem}
Suppose a system of polynomials satisfies a three-term recurrence relation
of the form \thetag{1.3} with $a_{k+1} > 0$ and $b_k \in {\Bbb R}$
$(k=0,1,2,\ldots)$ and initial conditions $r_{-1}(z) = 0$ and $r_0(z) = 1$, then
these polynomials are orthonormal in $L^2(\mu)$ for some positive
measure $\mu$ on the real line.
\endproclaim
This theorem is usually called Favard's theorem \cite{25} but it
is basically already in Stieltjes' memoir \cite{105, \S 11} for the case
of $J$-fractions obtained from contracting an $S$-fraction with positive
coefficients: he shows
that there is a positive linear functional $S$ such that $S(r_mr_n) = 0$
whenever $m \neq n$. The only thing that Stieltjes was missing was the Riesz
representation theorem which would enable one to express the linear functional
$S$ as a Stieltjes integral.

Hilbert's work on quadratic forms in infinitely many variables was much inspired
by Stieltjes' work on continued fractions \cite{43, p.\ 109}: {\it ``Die
Anwendungen
der Theorie  der quadratischen Formen mit unendlich vielen Variabeln sind nicht
auf die Integral\-glei\-chungen beschr\"ankt: es bietet sich nicht minder
eine Ber\"uhrung dieser Theorie mit der sch\"onen Theorie der Kettenbr\"uche von
Stieltjes ...''}. Stieltjes' theory is full of important ideas.
In Chapter V of \cite{105} Stieltjes gives a discussion on the
convergence of sums of the form $$    f_1(z)+f_2(z)+\cdots+f_n(z), $$
where $f_i(z)$ are analytic functions on the open unit disk $C_R$ with
center at the origin and radius $R$.  He proves a result which was later
also proved by Giuseppe Vitali in 1903 \cite{118}:
\proclaim{Theorem (Stieltjes-Vitali)}
Let $f_n$ be a sequence of analytic functions on a non\-empty connected open set
$\Omega$
of the complex plane. If $f_n$ is uniformly bounded on compact sets of $\Omega$
and if $f_n$ converges on a subset $E \subset \Omega$ that has an
accumulation point in $\Omega$, then $f_n$ converges uniformly on every compact
subset of $\Omega$.
\endproclaim
Paul Montel refers to this theorem as Stieltjes' theorem
\cite{67} and others refer to it as Vitali's theorem.
This result is very convenient in the study of
convergence of continued fractions
because quite often one is dealing with rational
fractions and one may be able to prove convergence on a set $E$ that is far
enough away from the poles of the rational fraction. The Stieltjes-Vitali
theorem then allows one to extend the asymptotic result to hold everywhere
except at the set containing all the poles.

The continued fraction \thetag{1.2} was later studied by Van Vleck \cite{116}
for $b_k$ arbitrary real numbers and $a_k^2$ arbitrary positive numbers. The
corresponding
measure is then not necessarily supported on $[0,\infty)$ and these continued
fractions are then closely related to Stieltjes integrals over
$(-\infty,\infty)$. The complete extension is due to Hamburger \cite{37}.
Van Vleck \cite{115} and Pringsheim \cite{83} \cite{84} have also given an
extension to complex coefficients.

For some good expositions on continued fraction we refer to the books
by Perron \cite{78}, Wall \cite{119}, Jones and Thron \cite{46} and
Lorentzen and Waadeland \cite{58}.

\head 2. Moment Problems \endhead

\subhead 2.1 The Stieltjes moment problem \endsubhead

In his fundamental work \cite{105, \S 24} Stieltjes introduced the following
problem: given an infinite sequence $\mu_k$ $(k=0,1,2,\ldots)$, find a
distribution of mass (a positive measure $\mu$) on the semi-infinite
interval $[0,\infty)$ such that
$$    \mu_k = \int_0^\infty x^k \, d\mu(x), \qquad k=0,1,2,\ldots . $$
Of course such a measure will not always exist for any sequence $\mu_k$ and
if such a measure exists, then it need not be unique. The Stieltjes moment
problem therefore has two parts
\roster
\item find necessary and/or sufficient conditions for the existence of a
solution of the moment problem on $[0,\infty)$,
\item find necessary and/or sufficient conditions for the uniqueness of the
solution of the moment problem on $[0,\infty)$.
\endroster
Chebyshev had previously investigated integrals and sums of the form
$$  \int_{-\infty}^{\infty} \frac{w(t)}{x-t} \, dt , \qquad \sum_{i=0}^{\infty}
\frac{w_i}{x-x_i}, $$
where $w(t)$ is a positive weight function and $w_i$ are positive weights.
Stieltjes integrals cover both cases and give a unified approach to the
theory. Chebyshev did not investigate a moment problem, but was interested
when a given sequence of moments determines the function $w(x)$ or the
weights $w_i$ uniquely. His work and the work of his student Markov
is very relevant, but Stieltjes apparently was unaware of it.
See Krein \cite{53} for some history related to the work of Chebyshev and
Markov.
Nevertheless Stieltjes' introduction of the moment problem is still
regarded as an important mathematical achievement.
The reason for the introduction of this moment problem is a close connection
between $S$-fractions or $J$-fractions and infinite series. If we make a
formal expansion of the function
$$   S(\mu;x) = \int \frac{d\mu(t)}{x+t},  $$
which is known as the Stieltjes transform of the measure $\mu$, then we find
$$  \int \frac{d\mu(t)}{x+t} \sim \sum_{k=0}^\infty (-1)^k
         \frac{\mu_k}{x^{k+1}} . $$
This series does not always converge and should be considered as an
asymptotic expansion.
On the other hand one can expand the function $S(\mu;x)$ also into a continued
fraction of the form \thetag{1.1} or \thetag{1.2}. The $n$-th approximant
of the $J$-fraction has the property that the first $2n$ terms in the expansion
$$    \frac{1}{a_1} \, \frac{p_{n-1}^{(1)}(x)}{p_n(x)}  \sim
    \sum_{k=0}^{\infty} (-1)^k \frac{m_k}{x^{k+1}}   $$
agree with those of the expansion of $S(\mu;x)$ i.e., $m_k=\mu_k$ for
$k=0,1,\ldots,2n-1$. This rational function is
therefore a (diagonal) Pad\'e approximant for $S(\mu;x)$.
If the infinite series is
given, then the  continued fraction is completely known whenever the measure
$\mu$ is known, provided the continued fraction converges.

Stieltjes gave necessary and sufficient conditions for the existence of a
solution of the Stieltjes moment problem:
\proclaim{Theorem}
If the Hankel determinants satisfy
$$  \vmatrix
    \mu_0 & \mu_1 & \cdots & \mu_n \\
    \mu_1 & \mu_2 & \cdots & \mu_{n+1} \\
    \vdots & \vdots & \cdots & \vdots \\
    \mu_n & \mu_{n+1} & \cdots & \mu_{2n}
    \endvmatrix      > 0, \qquad n \in {\Bbb N},  \tag 2.1 $$
and
$$  \vmatrix
    \mu_1 & \mu_2 & \cdots & \mu_{n+1} \\
    \mu_2 & \mu_3 & \cdots & \mu_{n+2} \\
    \vdots & \vdots & \cdots & \vdots \\
    \mu_{n+1} & \mu_{n+2} & \cdots & \mu_{2n+1}
    \endvmatrix      > 0, \qquad n \in {\Bbb N},  \tag 2.2 $$
then there exists a solution of the Stieltjes moment problem.
\endproclaim
If the moment problem has a unique solution then the moment problem is
determinate. If there exist at least two solutions then the moment problem is
indeterminate. Other terminology is also in use: determined/indetermined
and determined/undetermined.
Any convex combination of two solutions is another solution, hence
in case of an indeterminate moment problem there will always be an infinite
number of solutions. Stieltjes gave explicit examples of indeterminate moment
problems (see also Section 5.3) and he showed that a moment problem is
determinate if and only if the corresponding continued fraction \thetag{1.1}
converges for every $z$ in the complex plane, except for $z$ real and negative.
A necessary and sufficient
condition for a determinate moment problem is the divergence of the series
$\sum_{n=1}^{\infty} c_n$ where $c_n$ are the coefficients of the
$S$-fraction \thetag{1.1}. In case of an indeterminate moment problem Stieltjes
constructs two solutions as follows: let $P_n(z)/Q_n(z)$ be the
$n$-th convergent of the continued fraction \thetag{1.1}, then the limits
$$ \align
 \lim_{n \rightarrow \infty} P_{2n}(z) = p(z),
&\quad  \lim_{n \rightarrow \infty} P_{2n+1}(z) = p_1(z), \\
  \lim_{n \rightarrow \infty} Q_{2n}(z) = q(z),
&\quad  \lim_{n \rightarrow \infty} Q_{2n+1}(z) = q_1(z),
 \endalign   $$
exist, where $p,p_1,q,q_1$ are entire functions satisfying
$$    q(z)p_1(z)-q_1(z)p(z) = 1. $$
Stieltjes then shows that
$$    \frac{p(z)}{q(z)} = \sum_{k=1}^{\infty} \frac{r_k}{z+x_k} , \quad
   \frac{p_1(z)}{q_1(z)} = \frac{s_0}{z} + \sum_{k=1}^{\infty}
\frac{s_k}{z+y_k}. $$
The poles $x_k, y_k$ $(k=1,2,\ldots)$ are all
real and positive and the residues $r_k, s_k$ are all positive: this follows
because the zeros of the numerator polynomials interlace with the zeros of the
numerator polynomials and because all these zeros are real and negative.
These limits can thus be expressed as a Stieltjes integral
$$  \frac{p(z)}{q(z)} = \int_0^\infty \frac{d\mu(t)}{z+t}, \qquad
  \frac{p_1(z)}{q_1(z)} = \int_0^\infty \frac{d\mu_1(t)}{z+t}, $$
and both $\mu$ and $\mu_1$ are solutions of the moment problem with
remarkable extremal properties. This is one instance where it is clear
why Stieltjes  introduced the concept of a Stieltjes integral.

Not much work on the Stieltjes moment problem was done after Stieltjes'
death. One exception is
G. H. Hardy \cite{38} who considered the moments of a weight function $w(x)$ on
$[0,\infty)$ with restricted behaviour at infinity:
$$    \int_0^{\infty} w(x)e^{k\sqrt{x}} \, dx  < \infty, $$
for a positive value of $k$. He shows that the Stieltjes moment problem is
then always determinate and constructs the density from the series
$$  \sum_{n=0}^{\infty} \frac{\mu_n (-x)^n}{(2n)!} . $$
Hardy's proof avoids the use of continued fractions.

\subhead 2.2 Other moment problems \endsubhead

Nothing new happened until 1920 when Hamburger \cite{37} extended
Stieltjes'
moment problem by allowing the solution to be a measure on the whole real line
instead of the positive interval $[0,\infty)$. The extension seems
straightforward but the analysis is more complicated because the
coefficients of the continued fraction \thetag{1.1} may become negative or
vanish. Hamburger showed, using continued fraction techniques, that a
necessary
and sufficient condition for the existence of a solution of the Hamburger moment
problem is the positivity of the Hankel determinants \thetag{2.1}. He also shows
that a Hamburger moment problem may be indeterminate while the Stieltjes
moment problem with the same moments is determinate.

Nevanlinna \cite{70} introduced techniques of modern function
theory to investigate moment problems without using continued fractions.
His work is important because of the notion of extremal solutions, which were
first studied by him. M.\ Riesz \cite{85} \cite{86} gave a close connection
between the density of polynomials in $L^2$-spaces and moment problems:

\proclaim{Theorem}
Let $\mu$ be a positive measure on $(-\infty,\infty)$. If the Hamburger moment
problem for $\mu_k = \int x^k \, d\mu(x)$ $(k=0,1,2,\ldots)$ is determinate,
then polynomials are dense in $L^2(\mu)$. If the Hamburger
moment problem is indeterminate then the polynomials are dense in $L^2(\mu)$
if and only if $\mu$ is a Nevanlinna extremal measure.
\endproclaim
Berg and Thill \cite{9} have recently pointed out that this connection
is not any longer valid in higher dimensions by showing that there exist
rotation invariant measures $\mu$ on ${\Bbb R}^d, d > 1$ for which the
moment problem is determinate but for which polynomials are not dense in
$L^2(\mu)$.

In 1923 Hausdorff \cite{39} studied the moment problem for measures on a
finite interval $[a,b]$. The Hausdorff moment problem is always determinate
and conditions for the existence of a solution can be given in terms of
completely monotonic sequences.
The moment problem is closely related to quadratic forms of infinitely many
variables and operators in Hilbert space, as became clear from the work of
Carleman \cite{12} \cite{13} and Stone \cite{107}. Carleman established the
following sufficient condition for a determinate moment problem:
$$   \sum_{k=1}^\infty \mu_{2k}^{-1/2k} = \infty . $$
This is still the most general sufficient condition.
Karlin and his collaborators \cite{50} \cite{51} have approached the moment
problem through the geometry of convex sets and have shown that many
results can be interpreted in this geometrical setting.
Let me mention here that one can find excellent treatments of the moment problem
in the monograph of Shohat and Tamarkin \cite{92} and the book of Akhiezer
\cite{3}. Also of interest is the monograph by Krein and Nudelman
\cite{54}.

\subhead 2.3 Recent extensions of the moment problem \endsubhead

The most recent extension of the moment problem is to consider a doubly
infinite sequence $\mu_n$ $(n \in {\Bbb Z})$ and to find a positive measure
$\mu$ on $(-\infty,\infty)$ such that
$$  \mu_n = \int_{-\infty}^{\infty} x^n \, d\mu(x),  \qquad n \in {\Bbb Z}. $$
Such a moment problem is known as a strong moment problem. The strong
Stieltjes moment problem was posed and solved by Jones, Thron and Waadeland
in 1980 \cite{48} and again the solution is given in terms of the positivity
of certain Hankel determinants. These authors again use continued fractions,
but instead of the $S$- and $J$-fractions encountered by Stieltjes and
Hamburger, one  deals with another kind of fraction known as a $T$-fraction.
The strong Hamburger moment problem was handled by Jones, Thron and Nj\aa stad
in 1984
\cite{47}. Nj\aa stad \cite{71} gave another extension, known as the extended
moment
problem: given $p$ sequences $\mu_n^{(k)}$ $(n=1,2,\ldots; 1 \leq k \leq p)$ and
$p$ real numbers $a_1,a_2,\ldots,a_p$, does there exist a positive measure
$\mu$ on the real line such that
$$   \mu_n^{(k)} = \int_{-\infty}^{\infty} \frac{d\mu(t)}{(t-a_k)^n},
\qquad 1 \leq k \leq p, n \in \Bbb N\ ? $$
The solution is again given in terms of positive definiteness of a certain
functional.
Orthogonal polynomials play an important role in the Stieltjes and Hamburger
moment problem; for the strong moment problem a similar important role is played
by orthogonal Laurent polynomials and for the extended moment problem one deals
with orthogonal rational functions. The first place where orthogonal Laurent
polynomials are considered seems to be a paper by Pastro \cite{75}, where
an explicit example of the orthogonal Laurent polynomials with respect to the
Stieltjes-Wigert weight appears (see \S 5.3 for this weight).

\head 3. Electrostatic Interpretation of Zeros  \endhead

Stieltjes gave a very interesting interpretation of the zeros of Jacobi,
Laguerre and Hermite polynomials in terms of a problem of
electrostatic equilibrium. Suppose $n$ unit charges at  points
$x_1,x_2,\ldots,x_n$ are distributed in the (possibly infinite) interval
$(a,b)$.
The expression
$$    D(x_1,x_2,\ldots,x_n) = \prod_{1 \leq i < j \leq n} |x_i-x_j|  $$
is known as the discriminant of $x_1,x_2,\ldots,x_n$. If the charges
repell each other according to the law of logaritmic potential, then
$$    -\log D(x_1,x_2,\ldots,x_n) = \sum_{1 \leq i < j \leq n} \log
\frac{1}{|x_i-x_j|}   $$
is the energy of the system of electrostatic charges and the minimum of
this
expression gives the electrostatic equilibrium. The points $x_1,x_2,\ldots,x_n$
where the minimum is obtained are the places where the charges will settle
down. Stieltjes observed that these points are closely related to zeros of
classical orthogonal polynomials.

\subhead 3.1 Jacobi polynomials \endsubhead

Suppose the $n$ unit charges are distributed in $[-1,1]$ and that we add
two
extra charges at the endpoints, a charge $p>0$ at $+1$ and a charge $q>0$ at
$-1$.
Each of the unit charges interacts with the charges at $\pm 1$ and
therefore the electrostatic energy becomes
$$   L = -\log D_n(x_1,x_2,\ldots,x_n) + p \sum_{i=1}^n
\log \frac{1}{|1-x_i|} + q  \sum_{i=1}^n \log \frac{1}{|1+x_i|}  . \tag 3.1 $$
Stieltjes then proved the following result \cite{97} \cite{98} \cite{100}
\proclaim{Theorem}
The expression (3.1) becomes a minimum when $x_1,x_2,\ldots,x_n$ are the zeros
of the Jacobi polynomial $P_n^{(2p-1,2q-1)}(x)$.
\endproclaim

\demo{Proof:}
It is clear that for the minimum all the $x_i$ are distinct and different from
$\pm 1$.
For a minimum we need $\partial L/\partial x_k = 0$ $(1 \leq k \leq n)$
so that we have the system of equations
$$  \sum\Sb i=1 \\ i \neq k \endSb^n \frac{1}{x_i-x_k} - \frac{p}{x_k-1}
 - \frac{q}{x_k+1} = 0, \qquad 1 \leq k \leq n . $$
If we introduce the polynomial
$$   p_n(x) = (x-x_1)(x-x_2)\cdots(x-x_n), $$
then this is equivalent with
$$   \frac12 \frac{p_n''(x_k)}{p_n'(x_k)} + \frac{p}{x_k-1} + \frac{q}{x_k+1} =
0,  \qquad 1 \leq k \leq n.$$
This means that the polynomial
$$    (1-x^2)p_n''(x) + 2[q-p-(p+q)x] p_n'(x)  $$
vanishes at the points $x_k$ and since this polynomial is of degree $n$
it must be a multiple of $p_n(x)$. The factor is easily obtained by
equating the coefficient of $x^n$ and we have
$$    (1-x^2)p_n''(x) + 2[q-p-(p+q)x] p_n'(x) = -n[n+2(p+q)-1] p_n(x), $$
which is the differential equation for the Jacobi polynomial
$P_n^{(2p-1,2q-1)}(x)/c_n$, where $c_n$ is the leading coefficient of the Jacobi
polynomial. \qed
\enddemo
Stieltjes also found the minimum value. Hilbert \cite{42} also computed the
minimum value and Schur \cite{89} treated the case $p=q=0$ in detail.
Schur's
paper then led Fekete \cite{26} to define the transfinite diameter of a
compact set $K$ (with infinitely many points) in the complex plane. Take $n$
points $z_i \in K$ $(i=1,2,\ldots,n)$, and put
$$     d_n =  \max_{z_i \in K} D(z_1,\ldots,z_n)^{1/\binom{n}{2}} , $$
then $d_n$ is a decreasing and positive sequence \cite{26} \cite{111, Thm.\
III.21 on p.\ 71}. The limit of this sequence
is the transfinite diameter of $K$ and is an important quantity in logarithmic
potential theory (see \S 3.4). The transfinite diameter thus comes directly
from Stieltjes' work.

Consider the function
$$    \left( \prod_{i=1}^n |1-x_i|^{x-1}|1+x_i|^{y-1} \right)
     D_n(x_1,\ldots,x_n), $$
then Stieltjes' electrostatic interpretation gives the
$L_{[-1,1]^n}^\infty$-norm of this
function. The $L_{[-1,1]^n}^p$-norm of this function is also very famous and
is known as Selberg's beta integral \cite{90}. Actually Selberg
evaluated a multiple integral over $[0,1]^n$:
$$   \int_0^1 \ldots \int_0^1 D_n(t_1,\ldots,t_n)^{2z}
   \left( \prod_{i=1}^n t_i^{x-1}(1-t_i)^{y-1} \right) \, dt_1\ldots dt_n $$
$$  = \prod_{j=1}^n \frac{\Gamma(x+(j-1)z)\Gamma(y+(j-1)z)\Gamma(jz+1)}
     {\Gamma(x+y+(n+j-2)z)\Gamma(z+1)} , $$
but this integral can easily be transformed to an integral over
$[-1,1]^n$ which by an appropriate choice of the parameters $z,x,y$
becomes the desired $L_{[-1,1]^n}^p$-norm. This multiple integral
has many important applications e.g., in the statistical theory of high
energy levels (Mehta \cite{64}) but also in the algebraic theory of root
systems (Macdonald \cite{60}). Aomoto \cite{6} gave an elementary
evaluation of Selberg's integral and Gustafson \cite{35} computed some
$q$-extensions. Selberg's work was not inspired by Stieltjes, but it is directly
related to it.

\subhead 3.2 Laguerre and Hermite polynomials \endsubhead

A similar interpretation exists for the zeros of Laguerre and Hermite
polynomials.
Suppose the $n$ unit charges are distributed in $[0,\infty)$ and that we add
one
extra charge $p>0$ at the origin.
In order to prevent the charges from moving to $\infty$ we add the extra
condition that the centroid satisfies
$$   \frac1n \sum_{k=1}^n x_i \leq K ,  \tag 3.2  $$
with $K$ a positive number. The energy now is given by the expression
$$   L = - \log D_n(x_1,\ldots,x_n) + p \sum_{k=1}^n \log \frac1{x_k}. \tag
3.3 $$

\proclaim{Theorem}
The expression (3.3) together with the constraint (3.2) has a minimum
when $x_1,x_2,\ldots,x_n$ are the zeros of the Laguerre polynomial
$L_n^{(2p-1)}(c_nx)$, where $c_n=(n+2p-1)/K$. \endproclaim

If the $n$ unit charges are on $(-\infty,\infty)$ and if the moment of inertia
satisfies
$$            \frac1n \sum_{k=1}^n x_k^2 \leq L ,    \tag 3.4 $$
then

\proclaim{Theorem}
The expression $-\log D_n(x_1,x_2,\ldots,x_n)$ with constraint (3.4) becomes
minimal when $x_1,x_2,\ldots,x_n$ are the zeros of the Hermite polynomial
$H_n(d_nx)$, where $d_n=\sqrt{(n-1)/2L}$. \endproclaim

The proof of both statements is similar to the proof for the Jacobi case, except
that now we use a Lagrange multiplier to find the constrained minimum.
Mehta's book on Random Matrices \cite{64} gives an alternative way to
prove the results for Laguerre and Hermite polynomials.

In 1945 Siegel \cite{93} reproved the theorem for Laguerre polynomials and
applied it to improve the arithmetic-geometric mean inequality and to find
better bounds on algebraic integers. Siegels seems not to have been aware of
Stieltjes' work, but started from Schur's work \cite{89}.

\subhead 3.3 Extensions  \endsubhead

In \cite{99} Stieltjes generalizes this idea to polynomial solutions of the
differential equation
$$   A(x)y'' + 2B(x)y' + C(x)y = 0,    \tag 3.5  $$
where $A,B$ and $C$ are polynomials of degree respectively $p+1,p$ and
$p-1$. Such a differential equation
is known as a Lam\'e equation in algebraic form. In 1878 Heine asserted
that when $A$ and $B$ are given, there are in general exactly $\binom {n+p-1}n$
polynomials $C$ such that the differential equation has a solution which is a
polynomial of degree $n$. Stieltjes assumes that
$$ \frac{B(x)}{A(x)} = \sum_{k=0}^p \frac{r_k}{x-a_k},   $$
with $r_k > 0$ and $a_k \in {\Bbb R}$. One can then put charges $r_k$ at the
points $a_k$ and
$n$ unit charges at $n$ points $x_1,x_2,\ldots,x_n$ on the real line. Stieltjes
then shows that there are exactly $\binom {n+p-1}n$ positions
of electrostatic equilibrium, each  corresponding to one particular
distribution of the $n$ charges in the $p$ intervals
$[a_k,a_{k+1}]$ $(0 \leq k < p)$, and these charges are then at the points
$x_1,x_2,\ldots,x_n$ which are the $n$ zeros of the polynomial solution
of the differential equation. This result is now known as the Heine-Stieltjes
theorem \cite{110, Theorem 6.8 on p.\ 151}. The conditions imposed by Stieltjes
have been weakened by Van Vleck \cite{114} and P\'olya \cite{79}.
P\'olya allowed the zeros
of $A$ to be complex and showed that the zeros of the polynomial solution
of the differential equation will all belong to the convex hull of
$\{a_0,\ldots,a_p\}$. The location of the zeros  of the polynomial solution
is still under investigation now and interesting results and applications to
certain problems in physics and fluid mechanics are discussed in
\cite{4} \cite{5} \cite{123}.

Recently Forrester and Rogers \cite{27} and Hendriksen and van
Rossum \cite{40} have allowed the $n$ unit charges to
move into the complex plane. Forrester and Rogers consider a system of $2n$
particles of unit charge confined to a circle in the complex plane, say at the
points $e^{i\theta_j}$ and $e^{-i\theta_j}$ $(1 \leq j \leq n)$. At $\theta = 0$
(i.e., at the point $z=1$) a particle of charge $q$ is fixed and at $\theta =
\pi$ ($z=-1$) a particle of charge $p$. The energy of the system is now given by
$$  L = - q \sum_{k=1}^{2n} \log  |1-e^{i\theta_k}| - p \sum_{k=1}^{2n}
  \log |1+e^{i\theta_k}| - \sum_{1 \leq k < j \leq 2n} \log
|e^{i\theta_k}-e^{i\theta_j}| ,   \tag 3.6 $$
where
$$  0 < \theta_j < \pi, \quad \theta_j+\theta_{n+j} = 2\pi ,
    \qquad 1 \leq j \leq n .   \tag 3.7 $$
\proclaim{Theorem (Forrester and Rogers)}
The minimum of $L$ given in \thetag{3.6} subject to the constraints \thetag{3.7}
occurs when $\theta_j$ are the zeros of the trigonometric
Jacobi polynomial $P_n^{(q-\frac12,p-\frac12)}(\cos \theta)$.
\endproclaim
Forrester and Rogers also consider crystal lattice structures in which $n2^m$
particles of unit charge and $2^m$ particles of charge $q$ are distributed on
the unit circle, with one of the $q$ charges fixed at $\theta=0$. If one
requires that between every two $q$ charges there are $n$ unit charges
then the equilibrium position of the $n2^m$ particles of unit charge
occurs at the zeros of the Jacobi polynomial
$P_{n/2}^{(\frac{q-1}{2},-\frac12)}(\cos 2^m\theta)$ when $n$ is even and at the
zeros of $P_{(n-1)/2}^{(\frac{q-1}{2},\frac12)}(\cos 2^m\theta)$ when $n$ is
odd. The equilibrium position of the $2^m-1$ particles of charge $q$ occurs at
$\theta_k = \frac{2\pi k}{2^m}$ $(1 \leq k < 2^m)$.

Hendriksen and van Rossum \cite{40} have considered situations where other
special
polynomials come into play. Suppose $a > 0$ and that there is a charge $(a+1)/2$
at the origin and a charge $(c-a)/2$ at the point $1/a$. If $a \rightarrow
\infty$ one obtains a generalized dipole at the origin. Suppose now that
there are $n$
unit charges at points $z_1,z_2,\ldots,z_n$ in the complex plane, then the
electrostatic equilibrium in this generalized dipole field is obtained when
$z_1,\ldots,z_n$ are the zeros of the Bessel polynomial ${}_2F_0(-n,c+n;x)$.
Similar results can be obtained on so-called $m$-stars
$$   S_m = \{ z \in {\Bbb C} : z = \rho e^{\frac{2\pi k}{m}\,i},
  0 \leq \rho \leq r , k = 0,1,2,\ldots,m-1 \} . $$
Suppose that positive charges $q$ are placed at the endpoints $\rho=r$ of $S_m$
and a charge $p\geq 0$ is placed at the origin. If the points $z_1,\ldots,z_n$
$(n > m)$ in the complex plane all have a unit charge, then the
electrostatic equilibrium (assuming rotational symmetry) is obtained by choosing
$z_1,\ldots,z_n$ to be the zeros of the polynomial $f_n$ of degree $n$ that is a
solution of the differential equation
$$   (r^m-z^m)zy'' - 2[(p+qm)z^m-pr^m]y' = -n(n-1+2p+2qm)z^{m-1}y. $$
For particular choices of the parameters $p, r, m$ one then obtains well known
(orthogonal) polynomials.

\subhead 3.4 Logarithmic potential theory \endsubhead

Suppose that we normalize the electrostatic problem on $[-1,1]$ in such a way
that the total
charge is equal to 1. The $n$ charges then are equal to
$1/(n+p+q)$ and the charges at $1$ and $-1$ become respectively $p/(n+p+q)$ and
$q/(n+p+q)$. What happens if the number of particles $n$ increases? Clearly the
charges at the endpoints $\pm 1$ become negligible compared to the total charge
of the particles inside $[-1,1]$. This is the only place where $p$
and $q$ affect the distribution of the zeros, therefore it follows that the
asymptotic distribution of the charges in $(-1,1)$ i.e., the asymptotic
distribution of the zeros of Jacobi polynomials $P_n^{(2p-1,2q-1)}(x)$, is
independent of $p$ and $q$. By taking $p=q=1/4$ we deal with Chebyshev
polynomials of the first kind $T_n(x)$ with zeros $\{ \cos
\frac{(2j-1)\pi}{2n}: 1 \leq j \leq n\}$.
Let $N_n(a,b)$ be the number of zeros of $T_n(x)$
in $[a,b]$, then
$$ \align
   \frac{N_n(a,b)}{n} &= \sum_{a \leq \cos \frac{(2j-1)\pi}{2n} \leq b} \frac1n
\\ &= \int_{a \leq \cos t\pi \leq b} 1\, dt  + o(1) \\
                      &= \frac{1}{\pi} \int_a^b  \frac{dx}{\sqrt{1-x^2}} +
o(1). \endalign $$
Therefore the asymptotic distribution of the zeros of Jacobi polynomials is
given by the arcsin distribution and the relative number of zeros in $[a,b]$ is
$$    \frac1{\pi} \int_a^b \frac{dx}{\sqrt{1-x^2}} = \frac{1}{\pi} (\arcsin b -
\arcsin a ). $$
The surprising thing is that this is valid not only for Jacobi polynomials
but for a very large class of orthogonal polynomials on $[-1,1]$. The
arcsin distribution is actually an extremal measure in logarithmic potential
theory. Widom  \cite{120} \cite{121} and Ullman \cite{112} were probably the
first to connect logarithmic potential
theory and general orthogonal polynomials, even though some aspects such
as the transfinite diameter and conformal mappings had already appeared in
earlier work by Szeg\H{o} \cite{110, Chapter XVI}.
Let $K$ be a compact set in $\Bbb C$ and denote by $\Omega_K$ be the set
of all probability
measures on $K$. Define for $\mu \in \Omega_K$ the logarithmic energy by
$$  I(\mu) = \int_K \int_K \log \frac{1}{|x-y|} \,d\mu(x)\,d\mu(y), $$
then there exists a unique measure $\mu_K \in \Omega_K$ such that
$$   I(\mu_K) = \min_{\mu \in \Omega_K} I(\mu), $$
and this measure is the equilibrium measure (see e.g. \cite{111}).
When $K = [-1,1]$ then the equilibrium measure turns out to be the arcsin
measure. The capacity of the compact set $K$ is given by
$$   \text{cap}(K) = e^{-I(\mu_K)} , $$
and the capacity of a Borel set $B \in \Cal B$ is defined as
$$   \text{cap}(B) = \sup_{K \subset B,\ K \text{ compact}} \text{cap}(K), $$
(the capacity of $B$ is allowed to be $\infty$). Szeg\H{o} \cite{108}
showed that the capacity of a compact set $K$ is the same as the transfinite
diameter of this set, which we defined earlier.
The following result concerning the asymptotic
distribution of zeros of orthogonal polynomials is known (see e.g. \cite{94}):
\proclaim{Theorem}
Let $\mu$ be a probability measure on a compact set $K \subset \Bbb R$ such
that $$      \inf_{\mu(B)=1, B \in \Cal B} \text{cap}(B) = \text{cap}(K), $$
where $\Cal B$ are the Borel subsets of $K$, and suppose that $x_{k,n}$ $(1
\leq k \leq n)$ are the zeros of the orthogonal polynomial of degree $n$ for
the measure $\mu$. Then
$$   \lim_{n \rightarrow \infty} \frac1n \sum_{k=1}^n f(x_{k,n}) = \int_K f(t)
\, d\mu_K(t) $$
holds for every continuous function $f$ on $K$.
\endproclaim
When $K = [-1,1]$ then the conditions hold when $\mu$ is absolutely
continuous on $(-1,1)$ with $\mu'(x) > 0$ almost everywhere (in Lebesgue
sense). This includes all Jacobi weights. A very detailed account of
logarithmic potential theory and orthogonal polynomials can be found
in a forthcoming book by H. Stahl and V. Totik \cite{94}.

There is a similar generalization of the electrostatic interpretation of
the zeros of Laguerre and Hermite polynomials. This time we need to introduce
the energy of a measure in an external field $f$. If $K$ is a closed set in
the complex plane $\Bbb C$ and if the field $f: K \rightarrow [0,\infty)$ is
admissible i.e.,
\roster
\item $f$ is upper semi-continuous,
\item the set $\{ z \in K : f(z) > 0 \}$ has positive capacity ($\infty$ is
      allowed),
\item if $K$ is unbounded then $zf(z) \rightarrow 0$ as $|z|
\rightarrow \infty$ $(z \in K)$,
\endroster
then we define the energy integral in the field $f$ as
$$   I_f(\mu) = \int_K \int_K \log \frac{1}{|x-y|} \,d\mu(x)\,d\mu(y) - 2
 \int_K \log f(x) \, d\mu(x) . $$
Again there exists a unique measure $\mu_f$ such that
$$    I_f(\mu_f) = \min_{\mu \in \Omega_K} I_f(\mu) , $$
and this measure is the equilibrium measure in the external field $f$
\cite{34} \cite{65}. The
following result generalizes the electrostatic interpretation of the zeros of
Hermite polynomials (see e.g. \cite{34}):
\proclaim{Theorem}
Suppose that $x_{k,n}$ $(1 \leq k \leq n)$ are the zeros of the $n$-th degree
orthogonal polynomial with weight function $w(x)$ on $(-\infty,\infty)$. Suppose
that there exists a positive and increasing sequence $c_n$ such that
$$   \lim_{n \rightarrow \infty} w(c_nx)^{1/n} = f(x), \qquad x \in \Bbb R ,
                                \tag 3.8 $$
uniformly on every closed interval, with $f$ an admissible field, then
$$   \lim_{n \rightarrow \infty} \frac1n \sum_{k=1}^n g\left(
     \frac{x_{k,n}}{c_n} \right)  = \int g(x) \, d\mu_f(x)  , $$
for every bounded and continuous function $g$.
\endproclaim
 Again the asymptotic distribution of the (contracted) zeros of orthogonal
polynomials does not depend on the exact magnitude of the weight function $w$,
but only on the asymptotic behaviour given in (3.8). When $w(x) =
e^{-|x|^{\alpha}}$ --- the so-called Freud weights --- then one can take
$c_n = c(\alpha)n^{1/\alpha}$ with
$$     c(\alpha) = \left( \frac{\sqrt{\pi}\, \Gamma\left(
\frac{\alpha+1}{2} \right)}{\Gamma \left( \frac{\alpha}{2} \right)}
\right)^{1/\alpha},  $$
to find that $f(x) = e^{-|c(\alpha)x|^{\alpha}}$. The corresponding equilibrium
measure $\mu_f$ has support on $[-1,1]$ and is absolutely continuous with
weight function
$$     \mu_f'(t) = \frac1{\pi} \int_{|t|}^1 \frac{dy^{\alpha}}{\sqrt{y^2-t^2}},
\qquad -1 \leq t \leq 1 . $$
This is now known as the Nevai-Ullman weight.
Notice that the logarithm of the external field is the mathematical
counterpart of the constraints (3.2) and (3.4) for the Laguerre and Hermite
polynomials.

The fascinating aspects of logaritmic potential theory and zeros of orthogonal
polynomials are very much inspired by Stieltjes' observation that the zeros
of Jacobi, Laguerre and Hermite polynomials actually solve an equilibrium
problem in electrostatics.

\head 4. Markov-Stieltjes Inequalities \endhead

In his paper \cite{95} Stieltjes generalized the Gaussian quadrature formula,
which Gauss gave for the zeros of Legendre polynomials, to general weight
functions on an interval $[a,b]$. E. B. Christoffel had given this
generalization already seven years earlier \cite{19} \cite{30}, but
Stieltjes' paper is the first that makes a study of the convergence
of the quadrature formula. The Gaussian quadrature formula approximates
the integral
$$   \int_a^b \pi(x) \, d\mu(x),    \tag 4.1 $$
by appropriately summing $n$ function evaluations
$$  \sum_{j=1}^n \lambda_{j,n} \pi(x_{j,n}) . \tag 4.2 $$
This formula has maximal accuracy $2n-1$ i.e., the sum is equal to the integral
for all polynomials of degree at most $2n-1$, when the quadrature nodes
are the zeros $x_{j,n}$ $(1 \leq j \leq n)$ of the orthogonal polynomial
$p_n(x)$ of degree $n$ with orthogonality measure $\mu$, and the quadrature
weights $\lambda_{j,n}$ $(1 \leq j \leq n)$ are given by
$$   \lambda_{j,n} = \frac{-1}{a_{n+1}p_n'(x_{j,n})p_{n+1}(x_{j,n})} =
    \frac{1}{a_np_n'(x_{j,n})p_{n-1}(x_{j,n})} , $$
where we have used the recurrence relation \thetag{1.3}. These weights are known
as the {\it Christoffel numbers} and have important properties. One of the most
important properties is their positivity, which follows easily from
$$   \lambda_{j,n} =  \left\{ \sum_{k=0}^{n-1} p_k^2(x_{j,n}) \right\}^{-1}. $$
Stieltjes gives another remarkable property, namely
$$  \sum_{j=1}^{k-1} \lambda_{j,n} < \int_a^{x_{k,n}} d\mu(x) = \mu[a,x_{k,n})
\leq \mu[a,x_{k,n}] < \sum_{j=1}^k \lambda_{j,n} . \tag 4.3 $$
Stieltjes was unaware that Chebyshev had already conjectured these inequalities
in \cite{15} and that Chebyshev's student A. A. Markov had proved them  in
\cite{62}. Markov's paper appeared in 1884, the same year as Stieltjes' paper
\cite{95}, but in \cite{96} Stieltjes kindly acknowledges Markov to be the first
to prove the inequalities. He says however that his proof is independent of
Markov's proof since \cite{95} was submitted in May 1884 whereas Markov's paper
arrived at the library in September 1884. Szeg\H{o} \cite{110} gives three
proofs of \thetag{4.3}, combining the proofs of Stieltjes and Markov. The
inequalities \thetag{4.3} are nowadays known as the Markov-Stieltjes
inequalities. A related set of inequalities was proved
by K. Poss\'e \cite{80} \cite{81}. If $f: (a,b) \rightarrow {\Bbb R}$ is
such that $f^{(j)}(x) \geq 0$ for every $x \in (a,x_{k,n}]$
$(j=0,1,2,\ldots,2n-1)$ then
$$  \sum_{j=1}^{k-1} \lambda_{j,n}f(x_{j,n}) \leq
  \int_a^{x_{k,n}} f(x)\, d\mu(x) \leq  \sum_{j=1}^k \lambda_{j,n}f(x_{j,n}). $$
Stieltjes \cite{96} \cite{106} also gives other inequalities for the
Christoffel numbers e.g.,
$$   \sum_{j=1}^{k-1} \lambda_{j,n} < \sum_{j=1}^k \lambda_{j,n+1}
    < \sum_{j=1}^k \lambda_{j,n} .   \tag 4.4 $$

Stieltjes used the Markov-Stieltjes inequalities to show
that the sum \thetag{4.2} takes the form of
a Riemann-Stieltjes sum for the integral \thetag{4.1}, which makes Stieltjes the
first to prove convergence of Gaussian quadrature for continuous functions.
If $z \in {\Bbb C} \setminus [a,b]$ then the function $f(x) = 1/(z-x)$ is
continuous on $[a,b]$ and hence the Gaussian quadrature applied to $f$
converges to the Stieltjes transform of the orthogonality measure $\mu$. The
Gaussian quadrature formula for this function $f$ is a rational function
of $z$ and is exactly the $n$-th approximant ($n$-th convergent) of the
$J$-fraction for this Stieltjes transform and the convergence of the
Gaussian quadrature formula for $f$ therefore gives an important result
of Markov regarding the convergence of the diagonal
in the Pad\'e table for the Stieltjes transform of a positive measure.

Another important application of the Markov-Stieltjes inequalities is
a necessary and sufficient condition for determinacy of the moment problem: if
$$   \sum_{n=0}^\infty p_n^2(x) = \infty $$
for every real $x$ which is not a point of discontinuity of $\mu$, then
the moment problem for $\mu$ is determinate.
The Markov-Stieltjes inequalities are also very useful for estimations of the
rate of convergence of the Gaussian quadrature formula; the
Poss\'e-Markov-Stieltjes inequalities even give results for singular integrands
(Lubinsky and Rabinowitz \cite{59}). The estimation
of the distance between two succesive zeros of orthogonal (and also
quasi-orthogonal) polynomials can also be done using these inequalities.
From \thetag{4.3} one finds
$$    \lambda_{j,n} < \mu(x_{j+1,n},x_{j-1,n}) , \tag 4.5 $$
which allowed Stieltjes to deduce that $\lambda_{j,n}$ tends to zero when
$n \rightarrow \infty$ whenever the behaviour of $x_{j+1,n}-x_{j-1,n}$ is
known in terms of the measure $\mu$. Stieltjes gave the result for Legendre
polynomials.
Nevai \cite{68, p.\ 21} used the bound \thetag{4.5}
to show that for measures $\mu$ with compact support such that for every
$\epsilon > 0$ the set $\text{supp}(\mu) \cap (x-\epsilon,x+\epsilon)$
is an infinite set, there exists a sequence of integers $k_n$ $(n=1,2,\ldots)$
such that
$$      \lim_{n \rightarrow \infty} x_{k_n,n} = x, \quad
        \lim_{n \rightarrow \infty} \lambda_{k_n,n} = 0. $$
This shows that the zeros are dense in the derived set of $\text{supp}(\mu)$
and that the corresponding Christoffel numbers tend to zero.
For the isolated points in $\text{supp}(\mu)$ Nevai \cite{68, p.\ 156}
used \thetag{4.3} to show that
$$    \lim_{\epsilon \rightarrow 0+} \limsup_{n \rightarrow \infty}
    \sum_{|x-x_{k,n}|<\epsilon} \lambda_{k,n} = \mu(\{x\}), $$
for every $x \in {\Bbb R}$.
Freud \cite{28, p.\ 111} shows that for two
consecutive zeros in an interval $[c,d]$ for which
$$    0 < m < \frac{\mu(x,y]}{y-x} \leq M, \qquad x,y \in [c,d], $$
one has
$$     \frac{c_1}{n} \leq x_{j+1,n}-x_{j,n} \leq \frac{c_2}{n}, $$
where $c_1,c_2$ are positive constants. This is a slight extension of a result
by Erd\H{o}s and Tur\'an \cite{24}. Nevai \cite{68, p.\ 164} generalizes this
result by allowing $\mu'$ to have an algebraic singularity inside
$\text{supp}(\mu)$. If $\text{supp}(\mu)$ is compact, $\Delta \subset
\text{supp}(\mu)$, $t \in \Delta^o$ (the interior of the set $\Delta$)
and if $\mu$ is absolutely continuous
in $\Delta$ with
$$   c_1 |x-t|^\gamma \leq \mu'(x) \leq c_2 |x-t|^\gamma, \qquad \gamma > -1, $$
then
$$   \frac{c_3}{n} \leq  x_{k,n} - x_{k-1,n} \leq \frac{c_4}{n},  $$
whenever $x_{k,n} \in \Delta_1$ with $\Delta^1$ a closed subset of $\Delta^o$.
The Markov-Stieltjes inequalities are crucial to prove all these results.

\head 5. Special Polynomials \endhead
\subhead 5.1 Legendre polynomials \endsubhead

Stieltjes wrote a number of papers directly related to the Legendre polynomials
$P_n(x)$ for which
$$    \int_{-1}^1 P_n(x)P_m(x) \, dx = 0,  \qquad m \neq n. $$
He always uses the notation $X_n$ but here we will adopt the notation $P_n$
which is nowadays standard. In \cite{100} he uses the electrostatic
interpretation of the zeros of Jacobi polynomials to obtain monotonicity
properties of the zeros of Jacobi polynomials as a function of the parameters,
and from this one easily finds bounds for the zeros $x_{1,n} > x_{2,n} > \cdots
> x_{n,n}$ of the Legendre polynomials
$P_n(x) = P_n^{(0,0)}(x)$ in terms of the zeros of the Jacobi polynomials
$P_n^{(\frac12,-\frac12)}(x)$ and $P_n^{(-\frac12,\frac12)}(x)$ giving
$$      \cos \frac{2k\pi}{2n+1} < x_{k,n} < \cos \frac{(2k-1)\pi}{2n+1} ,
   \qquad 1 \leq k \leq n. $$
These bounds were already given by Bruns \cite{11} in 1881 and Stieltjes does
refer to Bruns' result, but Stieltjes goes on and shows that by using the zeros
of Chebyshev polynomials of the
first kind $T_n(x) = P_n^{(-\frac12,-\frac12)}(x)$ and
of the second kind $U_n(x) = P_n^{(\frac12,\frac12)}(x)$ one
may find better bounds:
$$    \cos \frac{k\pi}{n+1} < x_{k,n} < \cos \frac{(2k-1)\pi}{2n}, \qquad
1 \leq k \leq n/2. $$
A great deal of work has been done to obtain sharp bounds for zeros of
orthogonal
polynomials. The monotonicity of zeros of orthogonal polynomials depending on a
parameter is often used. Markov \cite{61} \cite{110, Theorem 6.12.1} gave a very
nice result concerning
the dependence of the zeros on a parameter $t$ which appears in the weight
function $w(x) = w(x;t)$. Two other methods for obtaining bounds for zeros of
orthogonal polynomials are the Sturm comparison theorem \cite{110, \S 6.3}
\cite{57} for solutions of Sturm-Liouville differential equations
and the Hellman-Feynman theorem \cite{45} of quantum chemistry.
See also \cite{44} for results on the monotonicity of zeros of orthogonal
polynomials.

Stieltjes made a very important contribution concerning the asymptotic
behaviour of Legendre polynomials. In 1878 Darboux \cite{23} gave an
asymptotic series for the Legendre polynomial:
$$   \align
    P_n(\cos \theta) = &2 a_n \sum_{k=0}^{m-1} a_k
\frac{1\cdot3\cdots(2k-1)}{(2n-1)(2n-3)\cdots(2n-2k+1)}  \\
   & \frac{\cos\left[ (n-k+\frac12)\theta - (k+\frac12)\pi/2 \right]}{(2\sin
\theta)^{k+\frac12}} + O(n^{-m-\frac12}), \qquad 0 < \theta < \pi,
  \endalign $$
which generalizes an asymptotic formula given by Laplace (when $m=1$). Here
$$    a_0 = 1, \quad a_k = \frac{1\cdot3\cdots(2k-1)}{2^k k!} . $$
The problem with this formula is that there is no closed expression or a bound
on the error term. Moreover the infinite series actually converges in the
ordinary sense when $\frac{\pi}{6} < \theta < \frac{5\pi}{6}$, but it converges
to $2P_n(\cos \theta)$ rather than $P_n(\cos \theta)$ (this ``paradox'' was
first pointed out by Olver \cite{72}). This is probably
the easiest example showing that asymptotic expansions need not converge
to the function that they approximate. The reason why things go wrong here is
that the formula is obtained by the so-called method of Darboux which consists
of obtaining asymptotic results of a sequence by carefully examining the
singularities on the circle of convergence of the generating function. The
generating function of Legendre polynomials has two singularities on the circle
of convergence, and at each singularity one picks up information on $P_n(\cos
\theta)$. This is probably the reason why the convergence of the infinite series
is to $2P_n(\cos \theta)$ rather than $P_n(\cos \theta)$. Stieltjes'
generalization of Laplace's asymptotic formula for the Legendre polynomials
does not suffer from either problem. Stieltjes' asymptotic expansion is
\cite{102} \cite{103}
$$  \align
   P_n(\cos \theta) = &\frac{4}{\pi} \frac{2^n n!}{3\cdot5\cdots(2n+1)} \\
        & \sum_{k=0}^{m-1} b_k \frac{\cos\left[(n+k+\frac12)\theta -
(k+\frac12)\pi/2\right]}{(2 \sin \theta)^{k+\frac12}} + R_m(\theta), \quad
 0 < \theta < \pi,
   \endalign $$
where
$$     b_0 = 1, \quad b_k = \frac{1^2\cdot3^2\cdots(2k-1)^2}{2^k k!
(2n+3)(2n+5)\cdots(2n+2k+1)}, $$
and the error $R_m(\theta)$ is bounded by
$$   |R_m(\theta)| < b_m \frac{4}{\pi} \frac{2^n n!}{3\cdot5\cdots(2n+1)}
 \frac{M}{(2 \sin \theta)^{m+\frac12}}, $$
where
$$ M = \cases 1/\cos \theta, & \text{if } \sin^2 \theta \leq \frac12, \\
    2 \sin \theta, & \text{if } \sin^2 \theta \geq \frac12. \endcases $$
This asymptotic expansion converges in the ordinary sense when $\frac{\pi}6 <
\theta < \frac{5\pi}{6}$ and it converges to $P_n(\cos \theta)$. Combined
with Mehler's asymptotic formula
$$  \lim_{n \rightarrow \infty} P_n(\cos \frac{\theta}{n} ) = J_0(\theta), $$
one then finds an asymptotic series for the Bessel function $J_0$ already
obtained by Poisson, but now with a bound on the error. Stieltjes also uses the
asymptotic series to obtain approximations of the zeros of the Legendre
polynomials.
The asymptotic theory of orthogonal polynomials (in particular classical
orthogonal polynomials) is very well developped nowadays, at least for
orthogonality on a finite interval. Szeg\H{o} has a very nice chapter on the
asymptotic properties of the classical polynomials \cite{110, Chapter VIII} and
that book is still a very good source for asymptotic formulas for Jacobi,
Laguerre and Hermite polynomials.

A third contribution of Stieltjes involving Legendre polynomials is his work on
Legendre functions of the second kind \cite{104}. The Legendre function of the
second kind can be defined by
$$   Q_n(x) = \frac12 \int_{-1}^1 \frac{P_n(y)}{x-y} \, dy , \qquad x \in {\Bbb
C} \setminus [-1,1], \tag 5.1 $$
so that
$$ Q_n(x) = \frac12 P_n(x) \log \left( \frac{x+1}{x-1} \right) -
P_{n-1}^{(1)}(x), \tag 5.2 $$
where $P_{n-1}^{(1)}(x)$ is the associated Legendre polynomial of degree $n-1$.
The integral representation cannot be used to define $Q_n(x)$ for $x \in
[-1,1]$ but by taking the appropriate limit, or the appropriate branch of the
logarithm, one can use \thetag{5.2} to define $Q_n(x)$ for $-1 < x < 1$.
Hermite \cite{41} had studied the zeros of $Q_n$ on $[-1,1]$ by making some
changes
of variables. Stieltjes works directly with $Q_n(x)$ as a function of the real
variable $x$ and shows that $Q_n(x)$ has $n+1$ zeros in $(-1,1)$ which
interlace with the zeros of the Legendre polynomial $P_n(x)$.
He also shows that there can be no zeros outside $[-1,1]$ by using a simple
property of Stieltjes transforms of positive weight functions.
Some of these results can easily be generalized to functions of the second kind
corresponding to general orthogonal polynomials \cite{113}.

\subhead 5.2 Stieltjes polynomials \endsubhead

In his last letter to Hermite \cite{8, vol.\ II, pp.\
439--441} Stieltjes considers the Legendre functions of the second kind
\thetag{5.1} and observes that
$$ \frac{1}{Q_n(z)} = E_{n+1}(z) + \frac{a_1}{z} + \frac{a_2}{z^2} + \cdots, $$
where $E_{n+1}(z)$ is a polynomial of degree $n+1$. This polynomial is now known
as the Stieltjes polynomial\footnote{The polynomial solutions of a Lam\'e
differential equation \thetag{3.5} are also known as Stieltjes polynomials but
we will not use that terminology.}.
Stieltjes gives the remarkable property
$$    \int_{-1}^1 P_n(x)E_{n+1}(x) x^k \, dx = 0, \qquad 0 \leq k \leq n, $$
which essentially means that $E_{n+1}(x)$ is orthogonal to all polynomials of
degree less than or equal to $n$ with respect to the oscillating weight
function $P_n(x)$ on $[-1,1]$. One may now wonder which
properties of ordinary orthogonal polynomials are still valid for $E_{n+1}(x)$
and Stieltjes conjectures that the zeros of $E_{n+1}(x)$ are real, simple and
belong to $[-1,1]$ and that they interlace with the zeros of $P_n(x)$. These
conjectures were later proved by Szeg\H{o} \cite{109}. Szeg\H{o} also extended
the idea to ultraspherical weights by considering the functions of the second
kind
$$  q_n^{\mu}(z) = \frac12 \frac{\Gamma(2\mu)}{\Gamma(\mu+\frac12)}
  \int_{-1}^1 (1-t^2)^{\mu-\frac12} \frac{P_n^{\mu}(x)}{z-x} \, dx, $$
where $P_n^{\mu}(x)$ is an ultraspherical polynomial of degree $n$.
One can then find
$$   \frac{1}{q_n^{\mu}(z)} = E_{n+1}^{\mu}(z) + \frac{a_1^\mu}{z} +
  \frac{a_2^\mu}{z^2} + \cdots, $$
where $E_{n+1}^{\mu}(z)$ is a polynomial of degree $n+1$. Szeg\H{o} shows
that
$$  \int_{-1}^1 (1-x^2)^{\mu-\frac12} P_n^\mu(x)E_{n+1}^\mu(x) x^k \, dx = 0,
 \qquad 0 \leq k \leq n, $$
thus generalizing the orthogonality of Stieltjes polynomials. The properties
of the zeros of $E_{n+1}^\mu(x)$ depend on the value of the parameter $\mu$.
If $0 < \mu \leq 2$ then the zeros of $E_{n+1}^\mu(x)$ are in $[-1,1]$, they
are real and simple and they interlace with the zeros of $P_n^\mu(x)$.
When $\mu < 0$ then some of the zeros are outside $[-1,1]$ and Monegato
\cite{66} has made
some computations showing that for $\mu \geq 4.5$ there can be complex zeros,
depending on the degree $n$.
More precise numerical information for Gegenbauer weights as well as for Jacobi
weights has been obtained by Gautschi and Notaris \cite{32}.
The construction of Stieltjes and Szeg\H{o} can be generalized by considering a
positive measure $\mu$ on ${\Bbb R}$. Suppose that $p_n(x;\mu)$
$(n=0,1,2,\ldots)$
are the orthogonal polynomials with respect to the measure $\mu$, then the
functions of the second kind are
$$    q_n(z;\mu) = \int \frac{p_n(x;\mu)}{z-x} \, d\mu(x), $$
and these are defined for $z \in {\Bbb C} \setminus \text{supp}(\mu)$. Define
the (general) Stieltjes polynomial $E_{n+1}(z;\mu)$ by
$$    \frac{1}{q_n(z;\mu)} = E_{n+1}(z;\mu) + \frac{a_1(\mu)}{z} +
   \frac{a_2(\mu)}{z^2}  + \cdots , $$
then one always has
$$   \int_{-1}^1 p_n(x;\mu) E_{n+1}(x;\mu) x^k \, d\mu(x) = 0, \qquad 0 \leq k
\leq n . $$

These Stieltjes polynomials turn out to have some importance in constructing an
optimal pair $(A,B)$ of quadrature formulas. Suppose we start with a quadrature
formula $A$ with $n$ nodes and a quadrature formula $B$ with $m$ nodes
$(m > n)$. In order to compute the error of formula $A$
one often assumes that the difference of the
results obtained by using $A$ and $B$ is proportional to the actual error of the
quadrature formula $A$. This means that one needs $n+m$ function evaluations to
compute the error of $A$. This implies that one has done $m$ extra function
evaluations which are not used in the evaluation of $A$ itself. Kronrod
\cite{55} suggested to extend formula $B$ to a formula with $n+m$ nodes in
such a way that the accuracy of $B$ is as high as possible.
For the Legendre weight on $[-1,1]$ one will find an optimal pair $(A,B)$ by
taking for $A$ the Gaussian quadrature with nodes equal to the zeros of the
Legendre polynomial $P_n(x)$ and for $B$ a quadrature formula with $2n+1$ nodes
at the zeros of $P_n(x)$ and the zeros of $E_{n+1}(x)$. The quadrature formula
$B$ then turns out to give a correct result for all polynomials of degree less
than or equal to $3n+1$ \cite{55, Theorem 6}.

In 1930 Geronimus \cite{33} slightly changes  Stieltjes' idea and considers the
Jacobi functions of the second kind
$$   Q_n^{(\alpha,\beta)}(z) = \int_{-1}^1 \frac{P_n^{(\alpha,\beta)}(x)}{z-x}
  (1-x)^{\alpha}(1+x)^{\beta} \, dx,  $$
with $P_n^{(\alpha,\beta)}(x)$ the Jacobi polynomial of degree $n$.
Geronimus observes that
$$  \frac{1}{Q_n(z) \sqrt{z^2-1}} = S_n(z) + \frac{c_1}{z} + \frac{c_2}{z}
    + \cdots , $$
with $S_n(z)$ a polynomial of degree $n$. Notice the extra factor $\sqrt{z^2-1}$
in the denominator on the left hand side. These polynomials satisfy the
remarkable property
$$    \int_{-1}^1 (1-x)^{\alpha}(1+x)^{\beta} P_n^{(\alpha,\beta)}(x)
S_n(x) T_k(x) \, dx = 0 ,   \qquad 0 < k \leq n, $$
and
$$    \int_{-1}^1 (1-x)^{\alpha}(1+x)^{\beta} P_n^{(\alpha,\beta)}(x)
S_n(x)  \, dx = 1. $$
Here $T_k(x)$ is the Chebyshev polynomial of the first kind of degree $k$.
Geronimus polynomials can be generalized to other weights on $[-1,1]$. The
interval $[-1,1]$ is important because it accounts for the factor $\sqrt{z^2-1}$
in the definition of the Geronimus polynomials. There is a relation between the
Geronimus polynomials $S_n(x)$ and
the Stieltjes polynomials $E_{n+1}(x)$ if one works with a weight function
on $[-1,1]$: if
$$   E_{n+1}(x) = \sum_{k=0}^{n+1}\!' \,c_{k,n} T_k(x) , $$
(the prime means to divide the first term by two)
is the  expansion of $E_{n+1}(x)$ in Chebyshev polynomials of the first kind,
then
$$   S_n(x) = \sum_{k=0}^n c_{k+1,n}U_k(x) $$
is the expansion of $S_n(x)$ in Chebyshev polynomials of the second kind.

Stieltjes and Geronimus polynomials and the related Gauss-Kronrod quadrature
are still being studied and we refer to Gautschi \cite{31},
Monegato \cite{66}, Peherstorfer
\cite{76} \cite{77} and Pr\'evost \cite{82} for more information.

\subhead 5.3 Stieltjes-Wigert polynomials \endsubhead

In his memoir \cite{105, \S 56} Stieltjes explicitly gives an example of a
moment problem on $[0,\infty)$ which is indeterminate. He shows that
$$ \int_0^\infty u^k u^{-\log u} \left[ 1 + \lambda \sin(2\pi \log u) \right]
 \, du = \sqrt{\pi} e^{\frac{(k+1)^2}{4}}    $$
is independent of $\lambda$ and therefore the weight functions
$$   w_\lambda(u) = u^{-\log u} \left[ 1 + \lambda \sin(2\pi \log u) \right],
        \qquad -1 \leq \lambda \leq 1 $$
all have the same moments which implies that this moment problem is
indeterminate.
Stieltjes gives the coefficients of the continued fraction \thetag{1.1}
$$  c_{2n} = (q;q)_{n-1} q^n, \quad c_{2n+1} =
\frac{q^{\frac{2n+1}{2}}}{(q;q)_n}  , $$
where $q=e^{-1/2}$ and
$$  (a;q)_0 = 1,\quad  (a;q)_n = (1-a)(1-aq)(1-aq^2)\cdots(1-aq^{n-1}) . $$
Both the series $\sum c_{2n}$ and $\sum c_{2n+1}$ converge since $0 < q < 1$,
which agrees with the theory worked out by Stieltjes. Later Wigert \cite{122}
extended this by considering the weight functions
$$   w_k(x) = e^{-k^2\log^2 x}, \qquad 0 < x <\infty, $$
which for $k=1$ reduce to the weight function considered by Stieltjes. If we set
$q=e^{-1/(2k^2)}$ then the orthogonal polynomials are given by
$$ p_n(x) = \sum_{j=0}^n \frac{(q^{-n};q)_j}{(q;q)_j} q^{j^2/2}(q^{n+1}x)^j, $$
and are known as Stieltjes-Wigert polynomials. The moment problem is
indeterminate whenever $0 < q < 1$, which means that there exist an infinite
number of measures on $[0,\infty)$ with the same moments. Askey \cite{7}
indicated that these polynomials are related to theta functions and shows
that the weight function
$$  w(x) = \frac{x^{-5/2}}{(-x;q)_{\infty}(-q/x;q)_{\infty}}, \qquad 0 < x <
\infty $$
has the same moments. This measure arises as a $q$-extension of the beta
density on $[0,\infty)$. Chihara \cite{16} \cite{18} has given many more
measures
which have the same moments as the weight function $w_k(x)$ given by Wigert.

The Stieltjes-Wigert polynomial $p_n(x)$ is a (terminating) basic hypergeometric
series. Such series are of the form $\sum c_j$ with $c_{j+1}/c_j$ a rational
function of $q^j$ for a fixed $q$ (for hypergeometric series this ratio is
a rational function of $j$). The first set of orthogonal polynomials which are
basic hypergeometric series was found by Markov in his thesis \cite{63}.
Except for a reference in Szeg\H{o}'s book \cite{110, \S 2.9}, this work was
overlooked and seems not to have led to any extensions. Markov's polynomials
are discrete extensions of Legendre polynomials and basic hypergeometric
extensions of discrete Chebyshev polynomials which are orthogonal on
$\{0,1,2,\ldots,N\}$ with respect to the uniform distribution. They are a
special case of polynomials considered by Hahn, which will be mentioned later.
The next basic hypergeometric orthogonal polynomials were introduced in 1894,
and there were two different examples that year. These are the
Stieltjes-Wigert polynomials (with $q=e^{-1/2}$) given by Stieltjes and
the continuous $q$-Hermite polynomials given by Rogers \cite{87}. Both
are basic hypergeometric extensions of Hermite polynomials but of a completely
different nature. Those of Rogers are orthogonal on $[-1,1]$ with respect
to the weight function
$$   w(x) = \frac{ \prod_{k=0}^\infty \left[ 1-2(2x^2-1)q^k+q^{2k}
\right]}{\sqrt{1-x^2}} . $$
A number of other examples were found before Hahn \cite{36} considered the
following problem: find all sets of orthogonal polynomials $p_n(x)$
$(n=0,1,2,\ldots)$ such
that $$   r_n(x) = \frac{p_{n+1}(x)-p_{n+1}(qx)}{x}, \qquad n=0,1,2,\ldots $$
is again a set of orthogonal polynomials. Earlier it had been shown that if
$p_n(x)$ $(n=0,1,2,\ldots)$ are orthogonal and $p_{n+1}'(x)$ $(n=0,1,2,\ldots)$
are orthogonal, then $p_n(x)$ $(n=0,1,2,\ldots)$ are either Jacobi, Laguerre
or Hermite polynomials (after a possible change of scale). It is easy to see
that the Stieltjes-Wigert polynomials are in the Hahn class.
The continuous $q$-Hermite polynomials of Rogers are not in the Hahn class, but
their analogous difference operator is a divided difference operator.
Basic hypergeometric series and orthogonal polynomials which are
terminating basic hypergeometric series are described in detail in the
book by Gasper and Rahman \cite{29}.
All of these polynomials arise in the study of quantum groups (see Koornwinder
\cite{52} and references there).

\subhead 5.4 Orthogonal polynomials related to elliptic functions \endsubhead

In Chapter XI of his memoir \cite{105} Stieltjes gives some examples of continued
fractions and
the corresponding moment problem. These examples (except for one) had already
been worked out in one of his previous papers \cite{101}. The continued fractions
are for the functions
$$  \align
   F_1(z,k) = \int_0^{\infty} \text{cn}(u,k)e^{-zu} \, du , &\qquad
   F_2(z,k) = \int_0^{\infty} \text{dn}(u,k)e^{-zu} \, du , \\
   F_3(z,k) = \int_0^{\infty} \text{sn}(u,k)e^{-zu} \, du , &\qquad
   F_4(z,k) = z \int_0^{\infty} \text{sn}^2(u,k)e^{-zu} \, du ,
    \endalign  $$
which are all Laplace transforms of the Jacobian elliptic functions
given by
$$  \align
    \text{cn}(u,k) &= \cos \varphi = \frac{2\pi}{kK} \sum_{n=1}^\infty
\frac{q^{n-\frac12}}{1+q^{2n-1}} \cos \frac{(2n-1)\pi u}{2K} ,  \\
    \text{sn}(u,k) &= \sin \varphi = \frac{2\pi}{kK} \sum_{n=1}^\infty
\frac{q^{n-\frac12}}{1-q^{2n-1}} \sin \frac{(2n-1)\pi u}{2K} ,  \\
    \text{dn}(u,k) &= \sqrt{1-k^2\sin^2 \varphi} = \frac{\pi}{2K}
+ \frac{2\pi}{K} \sum_{n=1}^\infty \frac{q^n}{1+q^{2n}} \cos \frac{n\pi u}{K} ,
   \endalign $$
with
$$    u = \int_0^\varphi \frac{d\theta}{\sqrt{1-k^2\sin^2 \theta}}, $$
and
$$   q = e^{-\pi K'/K}, \quad K(k) = \int_0^1
\frac{dx}{\sqrt{(1-x^2)(1-k^2x^2)}},
  \quad K'(k) = K(1-k^2). $$
The Chudnovsky's \cite{20, p.\ 197} pointed out that
these continued fractions are some of the very rare cases where both the
function and its continued fraction expansion are known explicitly. There
are quite a few cases known when the function is given in terms of
(basic) hypergeometric series and the numerators and denominators of the
convergents of the continued fraction are classical orthogonal polynomials
(in Askey's definition). The three-term recurrence relation then gives the
coefficients of the $J$-fraction. The functions $F_i(z,k)$ $(i=1,2,3,4)$ however
are not of (basic) hypergeometric type and the corresponding orthogonal
polynomials are therefore not classical. Nevertheless Stieltjes succeeded in
finding the continued fractions: he obtained $S$-fractions for $F_1$ and $F_2$
and $J$-fractions for $F_3$ and $F_4$.
His method consists of decomposing a quadratic form with infinitely many
variables as a sum of squares:
$$ \align
  \sum_{i=0}^\infty \sum_{j=0}^\infty a_{i+j} x_ix_j
  = & c_0 (x_0 + a_{0,1}x_1+a_{0,2}x_2 + \cdots)^2    \\
    &+ c_1(x_1 + a_{1,2}x_2+a_{1,3}x_3 + \cdots)^2   \\
    &+ c_2(x_2 + a_{2,3}x_3+a_{2,4}x_4 + \cdots)^2  + \cdots.
   \endalign  $$
The coefficients of the $J$-fraction of
$$   \sum_{n=0}^\infty \frac{(-1)^na_n}{z^{n+1}}  $$
are then determined by the coefficients $c_0,c_1,\ldots$ and $a_{i,i+1}$
$(i=0,1,2,\ldots)$. Such a decomposition can easily be made when
the function
$$   f(x) = \sum_{n=0}^\infty a_n \frac{x^n}{n!}  $$
satisfies an addition formula of the type
$$   f(x+y) = c_0f(x)f(y) + c_1 f_1(x)f_1(x) + c_2 f_2(x)f_2(y) + \cdots $$
where $f_m(x) = O(x^m)$, as Rogers \cite{88} pointed out when he was
reviewing Stieltjes' technique. The addition formulas for the Jacobian elliptic
functions then readily give the desired continued fractions.
The orthogonal polynomials that appear
are defined by the recurrence relations
$$    C_{n+1}(x) = xC_n(x) - \alpha_n C_{n-1}(x), \qquad
    D_{n+1}(x) = xD_n(x) - \beta_n D_{n-1}(x), $$
with
$$  \alpha_{2n} = (2n)^2k^2, \qquad \alpha_{2n+1} = (2n+1)^2, \qquad
    \beta_{2n} = (2n)^2, \qquad \beta_{2n+1} = (2n+1)^2k^2. $$
These polynomials have later been studied in detail by Carlitz \cite{14}.
The generating function for the
orthogonal polynomials satisfies a Lam\'e differential equation
$$   y'' + \frac12 \left\{ \frac1x + \frac1{x-1} + \frac1{x-a} \right\} y'
  + \frac{b-n(n+1)x}{4x(x-1)(x-a)} \, y = 0 , $$
with $n = 0$ (in general $n$ is an integer). There exist $2n+1$ values of the
parameter $b$ for which this Lam\'e equation has algebraic function solutions.
Stieltjes approach has been generalized to continued fraction expansions
for which the generating functions of the numerators and denominators of the
convergents satisfy a Lam\'e differential equation with $n \neq 0$. This has
been done by the Chudnovsky brothers \cite{20, pp.\ 197--201} \cite{21, \S 13}.
Some of these generalizations have interesting applications in number theory:
the irrationality and bounds on the measure of irrationality of some values of
complete elliptic integrals of the third kind can be obtained from these
continued fraction expansions.

\head Acknowledgements \endhead

I never quite realized how much work is needed to analyse Stieltjes' work a
century after his death. I have spent a lot of time in various libraries and
received a lot of help from the librarians. I would also like to thank various
colleagues for suggestions, comments and for pointing out omissions and
misinterpretations.
A sincere word of thanks in particular to Marcel de Bruin,
Ted Chihara, Walter Gautschi, Tom Koornwinder and Doron Lubinsky.
Of course nothing would
have been possible without the help of Gerrit van Dijk: many thanks for having
started this whole project.

\enddocument